\documentclass[letterpaper, 10 pt, conference]{ieeeconf}  

\IEEEoverridecommandlockouts                              

\usepackage{amsmath,amssymb,amsfonts}
\usepackage{graphicx,bm}
\usepackage[ruled,vlined]{algorithm2e}
\usepackage[colorlinks=true, bookmarks=true]{hyperref}
\AtBeginDocument{%
  \hypersetup{
    citecolor=red,
    linkcolor=blue,   
    urlcolor=blue}}

\allowdisplaybreaks

\title{\LARGE \bf
Optimal Trajectory Planning for Space Object Tracking \\ with Collision-Avoidance Constraints
}
\author{Saif R. Kazi$^{1}$, Harsha Nagarajan$^{1}$, Hassan Hijazi$^{1}$ and Przemek Wozniak$^{2}$
\thanks{*The authors gratefully acknowledge funding from LANL's Laboratory Directed Research \& Development program for the project ``20220392ER: Dynamic Coalition Management in Multi-Agent Sensor Networks''. \vspace{0.1cm}}
\thanks{$^{1}$ Applied Mathematics \& Plasma Physics (T-5), Los Alamos National Laboratory, Los Alamos, NM, USA.
        {\tt\small \{skazi, harsha, hlh\}@lanl.gov}}%
\thanks{$^{2}$ Space Remote Sensing \& Data Science (ISR-6), Los Alamos National Laboratory, Los Alamos, NM, USA.
        {\tt\small wozniak@lanl.gov}}%
}

\begin{document}

\maketitle
\thispagestyle{empty}
\pagestyle{empty}

\begin{abstract}

A control optimization approach is presented for a chaser spacecraft tasked with maintaining proximity to a target space object while avoiding collisions. The  target object trajectory is provided numerically to account for both passive debris and actively maneuvering spacecraft. Thrusting actions for the chaser object are modeled as discrete (on/off) variables to optimize resources (e.g., fuel) while satisfying spatial, dynamical, and collision-avoidance constraints. The nonlinear equation of motion is discretized directly using a fourth-order Runge-Kutta method without the need for linearized dynamics. The resulting mixed-integer nonlinear programming (MINLP) formulation is further enhanced with scaling techniques, valid constraints based on a perspective convex reformulation, and a combination of continuous relaxations of discrete actions with rounding heuristics to recover high-quality feasible solutions. This methodology enables efficient, collision-free trajectory planning over extended time horizons while reducing computational overhead. The effectiveness and practicality of the proposed approach is validated through a numerical case study.

\end{abstract}

\section{INTRODUCTION}
\label{sec:intro}

The problem of trajectory planning and optimization for space missions is as old as the space flight itself. A comprehensive review of the subject is presented in \cite{malyuta2021advances}. Continuing rapid advances in space exploration and new capabilities such as complex rendezvous and proximity operations (RPO) or modern spacecraft propulsion, call for new solutions to a multitude of challenging optimization problems that arise in this domain. At the same time, the emergence of powerful optimization methods combined with dramatic improvements in onboard computing is opening new possibilities for efficient and autonomous space flight control.

A major class of space vehicle control and trajectory optimization problems is concerned with one object (often referred to as the chaser) closely following another object (the target) in a (near)-Keplerian orbit around a central body, e.g. Earth. The two objects may potentially come into contact at the end (without crashing into any obstacles), as is the case for spacecraft docking or in situ debris mitigation. Exact discretization of the relative motion is possible with simplified linearized dynamics, such as the Hill–Clohessy–Wiltshire (HCW) equations \cite{HCW} for circular orbits or the linear time-varying Yamanaka–Ankerson state transition matrix \cite{YA} for elliptical orbits.
\begin{figure}[ht]
\centerline{\includegraphics[width=\columnwidth]{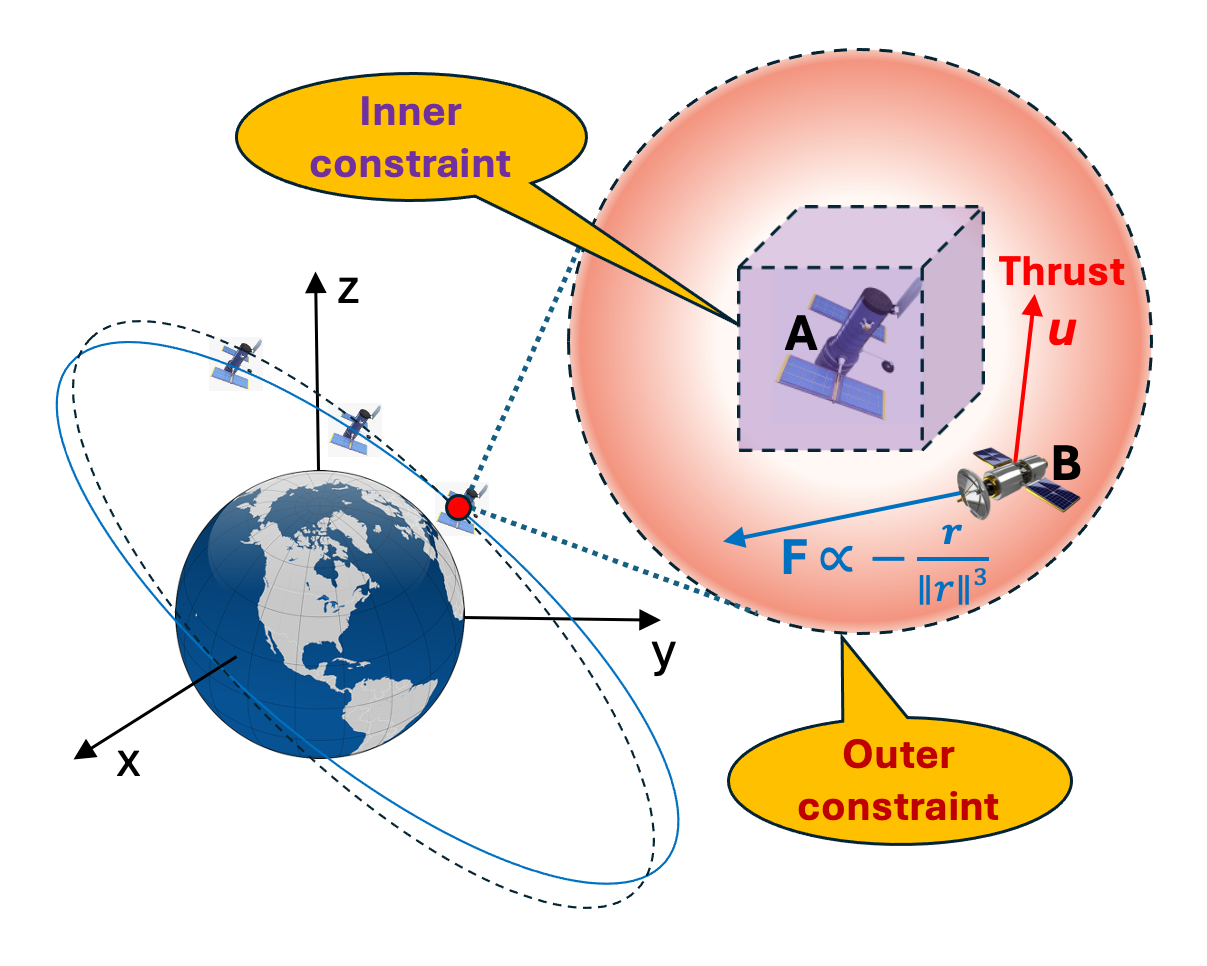}}
\caption{The constrained trajectory optimization problem involves the chaser spacecraft (B) maneuvering to closely follow the target object (A), whose motion is known, while avoiding collisions with the target and any other obstacles defined by ``keep out'' constraints.}
\label{fig:space_objects}
\end{figure}
Spatial constraints can be defined as polytope facets, ellipsoids, or a mix of both. Polytope exclusion regions lead to mixed-integer constraints, as first demonstrated in \cite{MIP_constraints} and \cite{MIP_example}. In the case of linearized dynamics, the resulting trajectory optimization can be solved using mixed-integer linear programming (MILP) methods \cite{MIP_constraints,MIP_example,tillerson2002co,basu2023computationally}.
In \cite{basu2023computationally}, the authors combine the target-assignment problem with the linearized dynamics of relative motion of satellites resulting in a large MILP model. Authors in \cite{tillerson2002co} implement their MILP solution in a nonlinear simulator with disturbances and actuator dynamics. An example rendezvous mission planning using mixed-integer nonlinear programming (MINLP) with linearized HCW equations is presented in \cite{ZHANG20111070}. Solving the exact nonlinear dynamics for satellites using collocation method along with collision avoidance constraints resulting in a nonlinear programming (NLP) was presented in \cite{lim2008trajectory}. It has also been shown that relative dynamics can be avoided in the inertial reference frame by lossless convexification of constraints and successively linearizing the dynamics \cite{Lu_and_Liu, Liu_and_Lu, Linearization_Filter}.

Artificial Intelligence and Machine Learning methods have been applied to a wide range of optimization problems in spacecraft guidance and control \cite{AI_and_ML}. Reinforcement Learning has emerged as a powerful unsupervised approach \cite{RL}. However, direct optimization approaches continue to have an advantage when it comes to providing performance guarantees such as the optimality gap and running time. In this study, we formulate the general fuel-minimization problem with exact nonlinear dynamics (up to the accuracy of the chosen discretization scheme) and discrete thruster actions as a large MINLP. We then propose reformulations using the perspective function to tighten the integer relaxation and decrease optimality gap. Unlike previous approaches, our method remains scalable to larger problem instances while ensuring feasible, near-optimal solutions.

Here, we consider constrained trajectory optimization for tracking space objects in a generalized case where the target motion is given numerically. Together with proximity constraints, this defines a series of approximate waypoints for the chaser trajectory to be optimized, with one waypoint at each time step (see Fig. \ref{fig:space_objects} for illustration). We formulate and solve the optimization problem in the inertial reference frame without the need for relative dynamics or linearization.

The nonlinear equation of motion is discretized using a fourth-order Runge-Kutta method. While this is a common numerical scheme for discretizing differential equations, it provides numerical stability and higher order accuracy of the optimal solution with respect to number of discretization steps. Other numerical discretization schemes can also be easily incorporated into our optimization model. For simplicity, the chaser spacecraft is assumed to be in a Keplerian orbit around the same central body as the target, this approach naturally allows the dynamics to be extended to include additional terms like departures from spherical symmetry of the gravitational field, residual atmospheric drag, and solar radiation pressure. Since the target trajectory appears in the proximity constraint, it may also reflect a variety of perturbations and non-gravitational forces, including active maneuvering.

Thruster action is quantized and modeled with on/off variables, leading to a Mixed-Integer Nonlinear Programming (MINLP) problem. We demonstrate that a feasible solution with a small and well understood optimality gap can be obtained in a relatively short time using modest computing hardware. The duration of the optimization interval is sufficient for potential onboard implementations even for relatively fast-moving objects in Low-Earth Orbits.

The main contributions of this paper are:
\begin{itemize}
    \item A MINLP model with nonlinear space dynamics for satellites with discrete thruster actions,
    \item A tight convex reformulation of the MINLP's continuous relaxation, obtained by applying established perspective-function techniques to reformulate the convex quadratic objective function,
    \item A simple and effective rounding strategy based on the convex relaxation solution, which yields feasible satellite control schedules with small optimality gaps.
\end{itemize}

The remainder of the paper is organized as follows: Section II presents the mathematical formulation of the optimization problem, followed by the solution algorithms in Section III. Section IV discusses the numerical results, and conclusions are provided in Section V.


\section{MATHEMATICAL FORMULATION}
\label{sec:form}

We derive the optimization model for the space object tracking by modeling the dynamics of the secondary object using Newton's law of gravity. For illustrative purposes, we present the system of equations in vector format and denote vector quantities in bold letters.

\subsection{Model Assumptions} 
Our optimization model assumes an Earth-Centered coordinate system as a non-accelerating inertial frame of reference, enabling a simplified two-body Keplerian dynamics model while neglecting Earth's revolution and higher-order gravitational harmonics. To isolate control-driven trajectory effects, we exclude external perturbations such as solar radiation pressure, third-body gravitational influences, and atmospheric drag. Both primary and secondary satellites are modeled as rigid bodies with static mass properties, disregarding structural flexure, fuel slosh, and mass reduction due to propellant consumption. Collision avoidance is enforced through deterministic geometric constraints, assuming ideal state knowledge of the system. These assumptions provide a tractable baseline for optimizing collision-free trajectories, ensuring that the thrust control strategy remains the primary variable influencing trajectory evolution.

\subsection{Orbital Dynamics Model}
The equation of motion for the secondary object (chaser) reads
\begin{align}\label{eq:newtons_law}
    \ddot{\bm{r}} = -\mu \frac{\bm{r}}{\|\bm{r}\|^3} + \bm{u},
\end{align}
where $\bm{r} = [x,y,z]$,  $\bm{u} = [u_x,u_y,u_z]$ and $\|\bm{r}\| = \sqrt{x^2 + y^2 + z^2}$. 

The equation of motion represented by Eq.~\eqref{eq:newtons_law} consists of two terms: The first term denotes Earth's gravitational pull experienced by the object using Newton's law of gravitation, where $\mu$ is the geocentric gravitational constant equal to $398600.436 \enskip \text{km}^3\text{s}^{-2}$. The second term represents control thrust, i.e. force per unit mass generated by the propulsion mechanism of the secondary object (e.g., ion-thrusters).

The radial distance from Earth's center $\bm{r}$ and thrust per unit mass $\bm{u}$ are vectors with individual components $x,y \enskip \text{and} \enskip z$. 

\subsection{Coupled First-order ODE}
Given $\bm{v} = [v_x,v_y,v_z]$,
\begin{subequations}\label{eq:first_order_ODE}
\begin{align}
    & \dot{\bm{r}}  = \bm{v}, \\
    & \ddot{\bm{r}}  = \dot{\bm{v}} = -\mu \frac{\bm{r}}{\|\bm{r}\|^3} + \bm{u},
\end{align}
\end{subequations}
we decompose the second-order differential equation of motion in \eqref{eq:newtons_law} into a system of coupled first-order differential equations, as shown in Eq.~\eqref{eq:first_order_ODE}. This approach facilitates numerical discretization and allows for the application of higher-order methods, improving accuracy. The velocity vector $\bm{v}$ represents the spacecraft's velocity relative to an inertial reference frame of reference centered at the Earth. Consequently, the system's state variables are the radial distance/position vector $\bm{r}$ and the velocity vector $\bm{v}$, while the thrust vector $\bm{u}$ represents the control variable.

\subsection{Discretization}

We use the fourth-order Runge Kutta (or RK4) method to discretize the coupled ordinary differential equations into a system of nonlinear equations \cite{ascher1998computer}. \\

\noindent \textbf{Runge-Kutta method}:
We define the time set $\mathbb{T} = \{ t_0, t_1,\ldots,t_f \}$ where the time horizon $[t_0, t_f]$ is uniformly divided into $N$ time intervals as  $t_{i+1} = t_i + \Delta t, i = 0,\ldots,N-1$. The state and control variables are discretized as $\mathcal{R} =  \{ \bm{r_0}, \bm{r_1},\ldots, \bm{r_f}\}$ and $\mathcal{V} =  \{ \bm{v_0}, \bm{v_1},\ldots, \bm{v_f}\}$, where the discrete values are approximation of the continuous variables at specific time instances i.e.  $\bm{r}(t_i) \approx \bm{r}_i$ and $\bm{v}(t_i) \approx \bm{v}_i$.
\begin{subequations}\label{eq:RK-discrete}
For $i = 0....N-1$
\begin{align}
    & \bm{r}_{i+1} = \bm{r}_i + \frac{\Delta t}{6}(\bm{k}_{1,i} + 2 \bm{k}_{2,i} + 2 \bm{k}_{3,i} + \bm{k}_{4,i}), \\
    & \bm{v}_{i+1} = \bm{v}_{i} + \frac{\Delta t}{6}(\bm{k}_{1,v,i} + 2 \bm{k}_{2,v,i} + 2 \bm{k}_{3,v,i} + \bm{k}_{4,v,i}),    
\end{align}   
\end{subequations}

The RK4 scheme is based upon a Taylor series with auxiliary variables ($\bm{k}$) as shown in \eqref{eq:RK-discrete} used to derive the discretized state equations. The Runge-Kutta variables for space $\bm{k}_{m,i} = [k_x,k_y,k_z]$ and velocity $\bm{k}_{m,v,i} = [k_{v_x},k_{v_y},k_{v_z}]$ are calculated as follows, where $m = 1,2,3,4$:
\begin{subequations}
\label{eq:RK-space}
\begin{align}
    & \bm{k}_{1,i} = \bm{f_r}(t_i,\bm{r}_i) = \bm{v}_i, \\
    & \bm{k}_{2,i} = \bm{f_r}(t_i + \Delta t/2, \bm{r}_i + k_1 \Delta t/2) = \bm{v}_i + \bm{k}_{1,v,i}\Delta t/2, \\ 
    & \bm{k}_{3,i} = \bm{f_r}(t_i + \Delta t/2, \bm{r}_i + k_2 \Delta t/2) = \bm{v}_i + \bm{k}_{2,v,i}\Delta t/2, \\
    & \bm{k}_{4,i} = \bm{f_r}(t_i + \Delta t, \bm{r}_i + k_3 \Delta t) = \bm{v}_i + \bm{k}_{3,v,i}\Delta t,
\end{align}  
where $\bm{f_r}(t,\bm{r}) = \dot{\bm{r}} = \bm{v}$,
\end{subequations}

\noindent Similarly for velocity-based Runge-Kutta variables
\begin{subequations}\label{eq:RK-velocity}
\begin{align}
    & \bm{k}_{1,v,i} = \bm{u}_i - \mu \bm{r}_i d_{1,i}, \\
    & \bm{k}_{2,v,i} = (\bm{u}_i + \bm{u}_{i+1})/2 - \mu (\bm{r}_i + \bm{k}_{1,i}\Delta t/2) d_{2,i}, \\ 
    & \bm{k}_{3,v,i} = (\bm{u}_i + \bm{u}_{i+1})/2 - \mu (\bm{r}_i + \bm{k}_{2,i}\Delta t/2) d_{3,i}, \\
    & \bm{k}_{4,v,i} = \bm{u}_{i+1} - \mu (\bm{r}_i + \bm{k}_{3,i}\Delta t) d_{4,i},
\end{align}    
here $\bm{f_v}(t,\bm{r}) = \bm{u} - \mu \bm{r} d$,
\end{subequations}

\noindent The inverse cubic distance variables ($d_1,d_2,d_3$ and $d_4$) are calculated as:
\begin{subequations}\label{eq:inverse cubic distance}
\begin{align}
    & d_{1,i} = (x_i^2 + y_i^2 + z_i^2)^{-3/2}, \\
    & d_{2,i} = ((x_i+k_{1,x,i}\Delta t/2)^2 + (y_i+k_{1,x,i}\Delta t/2)^2 + \\ 
    & \nonumber (z_i+k_{1,x,i}\Delta t/2)^2)^{-3/2}, \\
    & d_{3,i} = ((x_i+k_{2,x,i}\Delta t/2)^2 + (y_i+k_{2,x,i}\Delta t/2)^2 + \\ 
    & \nonumber (z_i+k_{2,x,i}\Delta t/2)^2)^{-3/2}, \\
    & d_{4,i} = ((x_i+k_{3,x,i}\Delta t)^2 + (y_i+k_{3,x,i}\Delta t)^2 + \\ 
    & \nonumber (z_i+k_{3,x,i}\Delta t)^2)^{-3/2}
\end{align}   
\end{subequations}

\subsection{Initial Conditions}

The initial conditions for the position and velocity of the secondary object are defined using the position and velocity of the primary object with some pre-specified distance ($\Delta$) between the two objects.
\begin{subequations}\label{eq:initial conditions}
\begin{align}
    & \bm{r}_0 = \bm{\Bar{r}}_0 + \bm{\Delta}, \\
    & \bm{v}_0 = \bm{\Bar{v}}_0,
\end{align}   
\end{subequations}
where $\bm{\bar r}$ and $\bm{\bar v}$ represent the primary object's initial location and velocity respectively. Note that the initial conditions are predefined and treated as fixed parameters in the optimal space object tracking problem.

\subsection{Thruster Constraints}

The secondary spacecraft thrusters are modeled as discrete control actions, meaning they can be turned ``on'' for limited time intervals, represented by binary variables $b_i \in \{0,1\}, \ \forall i = 0,\ldots,N$. 

\begin{itemize}
    \item When \( b_i = 1 \), the thruster can produce a nonzero control input, constrained within a predefined range:  
    $$\bm{u}_{\min} \leqslant \bm{u}_i \leqslant \bm{u}_{\max},$$
    \item When \( b_i = 0 \), the thruster remains off: $\bm{u}_i = 0 $. 
\end{itemize}

\begin{subequations}
\label{eq:thruster}
To enforce these conditions, we introduce an equivalent constraint that ensures the control input is only active when the corresponding binary variable is set to 1: 
\begin{align}\label{eq:big-M}
    b_i \bm{u}_{\min} \leqslant \bm{u}_i \leqslant b_i \bm{u}_{\max},
\end{align}

To constrain the number of times the thrusters can be activated within a given time period, we impose the following linear constraint:  
\begin{align}  
    \label{eq:thruster-budget}  
    \sum_{i=0}^{N} b_i \leqslant N_b  ,
\end{align}  
By restricting the total number of thruster activations to a predefined budget $N_b$, we effectively manage fuel consumption and operational efficiency. This constraint also reduces mechanical wear of thrusters, and applies thruster actions during critical maneuvers rather than dispersing it across multiple low-thrust intervals.
\end{subequations}

\subsection{Proximity Constraints}
Since the secondary spacecraft is tracking the primary object, proximity constraints are imposed to ensure the distance between the two objects remains within specified bounds, as modeled below:
\begin{align}\label{eq:proximity constraint}
    \Delta_{\min} \leqslant \| \bm{r - \bar r} \|_q \leqslant \Delta_{\max},
\end{align}
Depending on the norm type (\(q = 1\) or \(2\)), the constraint is reformulated into a sum of absolute terms or a quadratic constraint, respectively. The parameters $\Delta_{\min}$ and $\Delta_{\max}$ define the boundaries of the inner and outer regions shown in Fig.~\ref{fig:space_objects}.

\subsection{Objective Function}
The primary goal for the problem at hand is to achieve a specific orbit or proximity to a particular target object, with fuel consumption being the dominant factor. Thruster output is directly related to the fuel used. Therefore, we minimize the square of 2-norm of the thruster outputs as the objective function, as shown below:
\begin{align}\label{eq:objective}
    \min \left( \frac{1}{T} \int \|\bm{u}\|^2_2 \ dt\right) \ \approx \ \min \left( \frac{1}{N}  \sum_i^N \|\bm{u}_i\|^2 \right),
\end{align}
The time discretization enables the approximation of the continuous integral of the squared norm of the thrust vector using a Riemann sum, where the integral is replaced by a finite summation of the squared thruster values evaluated at discrete time points. This approximation converges to the exact integral as $N \to \infty$  under appropriate conditions on the discretization scheme. 

\subsection{Optimal Space Object Tracking Problem}

We formulate the optimal space object tracking problem as a MINLP with the discretized space dynamics, initial conditions, thruster dynamics and proximity condition as constraints. The space thruster outputs ($\bm{u}_i$) and the binary variables ($b_i$) are the degrees of freedom to the MINLP.
\begin{flalign*}
    \mathtt{MINLP} ~~ \triangleq ~~ & \min ~~ \text{Eq. \eqref{eq:objective}} ~~~ \text{ subject to: } \\
    & \text{Eq. \eqref{eq:RK-discrete}} \text{ -- (Runge-Kutta Discretized Equations)}, \\
    & \text{Eq. \eqref{eq:RK-space}, \eqref{eq:RK-velocity}, \eqref{eq:inverse cubic distance}} \text{ -- (Runge-Kutta variables)}, \\ 
    & \text{Eq. \eqref{eq:initial conditions}} \text{ -- (Initial Conditions)}, \\ 
    & \text{Eq. \eqref{eq:thruster}} \text{ -- (Thruster Constraints)}, \\ 
    & \text{Eq. \eqref{eq:proximity constraint}} \text{ -- (Proximity and Avoidance Constraints)},\\
    & \bm{r}_i, \bm{v}_i, \bm{u}_i \in \mathbb{R}^3, b_i \in B := \{0,1\}, \quad i = 0,\ldots,N
\end{flalign*}

\section{ALGORITHMS}
\label{sec:algorithm}

As presented in the previous section, the optimal space object tracking problem is modeled as a MINLP problem, which is $\mathcal{NP}$-hard in computational complexity. The complexity arises from the exponential increase in the number of feasible thrusting combinations for the binary variables as the number of time intervals ($N$) increases in conjunction with nonlinear, non-convex constraints. Additionally, the varying magnitudes of the variables—$\bm{r} \sim \mathcal{O}(10^3)$ km, $\bm{v} \sim \mathcal{O}(10^1)$ km/s, and $\bm{u} \sim \mathcal{O}(10^{-5})$ km\,s$^{-2}$ — necessitate that a good starting point is supplied to the solver. This is achieved by solving a simulation problem first.

\subsection{Simulation-based initialization}

The simulation problem is solved as an optimization problem with zero degrees of freedom. This is done by setting the values of the control variables equal to zero.
\begin{align}\label{zero control}
    \bm{u}_i = \bm{0},
\end{align}

We also remove variable bounds and inequality constraints resulting in a square system of nonlinear equations which is solved using a Newton-Raphson-based algorithm. 
\begin{flalign*}
    \mathtt{SIM} ~~ \triangleq ~~ &  \text{Eq. \eqref{eq:RK-discrete}, \eqref{eq:RK-space}, \eqref{eq:RK-velocity}, \eqref{eq:inverse cubic distance}, \eqref{eq:initial conditions}}, \qquad \bm{r}_i, \bm{v}_i \in \mathbb{R}^3
\end{flalign*}

We will denote the simulation solution as $\bm{x}_{\text{sim}} = \{\bm{r},\bm{v},\bm{u},\bm{k},\bm{k}_v\}$ and use it as the initial guess for solving the optimization problem with non-zero control variables.

\subsection{Continuous Relaxation of Binaries}
\label{subsec:bin_relax}
General MINLPs, as presented here, are very difficult to solve to optimality due to the nonlinear and nonconvex constraints. Even state-of-the-art global optimization solvers \cite{tawarmalani2005polyhedral,nagarajan2019adaptive} struggle to solve these problems within a reasonable optimality gap. To address this issue, we relax the binary variable \( b_i \in \{0,1\} \ \forall i \) to a continuous relaxation, as given by the constraint \( 0 \leqslant b_i \leqslant 1 \ \forall i=0,\ldots, N \).

The resulting continuous nonlinear program (NLP) after binary relaxation is still a challenging non-convex optimization to solve. However, we resort to solving it \textit{locally} using interior point methods. Specifically, we utilize Ipopt, which efficiently solves large-scale NLPs using a primal-dual interior-point method \cite{wachter2006implementation}. Ipopt leverages sparse matrix solvers and a filter line search method that balances feasibility and optimality. Additionally, Ipopt supports warm-starting solutions to an NLP, which can be effective in solving constrained problems of interest in this paper.
\begin{flalign*}
    \mathtt{RNLP} ~~ \triangleq ~~ & \min ~~ \text{Eq. \eqref{eq:objective}} ~~~ \text{ subject to: } \\
    & \text{Eq. \eqref{eq:RK-discrete}, \eqref{eq:RK-space}, \eqref{eq:RK-velocity}, \eqref{eq:inverse cubic distance}, \eqref{eq:initial conditions}, \eqref{eq:thruster}, \eqref{eq:proximity constraint}}, \\ 
    & \bm{r}_i, \bm{v}_i, \bm{u}_i \in \mathbb{R}^3, 0 \leqslant b_i \leqslant 1, \quad i = 0,\ldots,N
\end{flalign*}

\subsection{Perspective Convex Reformulation}
\label{sec:perspective}
Following the continuous relaxation from section \ref{subsec:bin_relax}, this procedure may yield fractional solutions for $b_i$, resulting in weak lower bounds for the MINLP. Moreover, these fractional thrusting values of $b_i$ can be infeasible to the trajectory optimization problem. 

A natural way to tighten the relaxation is to replace the nonconvex binary requirements with a semidefinite‐programming (SDP) relaxation \cite{axehill2010convex}. In practice, however, the resulting SDP-constrained continuous NLP is computationally formidable, and—crucially for the problem sizes studied here—no state-of-the-art solver can solve it to high accuracy within reasonable time. Instead, we adopt a tractable yet nearly as tight approach based on perspective reformulation. By applying the perspective mapping to the convex quadratic objective function in \eqref{eq:objective}, we obtain the following second-order-cone–representable relaxation that drives the solution toward the binary set $B:=\{0,1\}$:
\begin{subequations}
\begin{align}
    \min ~~ & \frac{1}{N} \sum_{i=0}^N \phi_i ~~ \text{subject to} \label{eq:perspective-cuta}\\ 
    & b_i \phi_i \geqslant \|\bm{u}_i\|^2, ~~ i = 0,\ldots,N. \label{eq:perspective-cutb}
\end{align}
\label{eq:perspective}
\end{subequations}
For every discretization $i$, the perspective reformulation in \eqref{eq:perspective} captures the tightest relaxation, which is the convex hull of the mixed-integer set \(S^i\) with three variables. Treating \(\|\bm{u}_i\|\) as a scalar real value, the set \(S^i = S^i_0 \cup S^i_1\), where
\begin{align}
S^i_0 = & \left\{ (\|\bm{u}_i\|, \phi_i, b_i) \in \mathbb{R}^2 \times B : \right. \nonumber \\  
& \left. \quad b_i = 0, \ \phi_i \geqslant 0, \  \|\bm{u}_i\| = 0 \right\}, \ \text{and} \\ \nonumber\\
S^i_1 = & \left\{ (\|\bm{u}_i\|, \phi_i, b_i) \in \mathbb{R}^2 \times B : \ b_i = 1, \ \phi_i \geqslant \|\bm{u}_i\|^2, \right. \nonumber \\  
& \quad \left. \|\bm{u}_{\min}\| \leqslant \|\bm{u}_i\| \leqslant  \|\bm{u}_{\max}\| \right\}.
\end{align}
By applying the perspective of the convex quadratic function in the set \(S^i\), the following convex constraint is obtained \cite{boyd2004convex}:
\begin{align}
\frac{\phi_i}{b_i} \geqslant \left( \frac{\|\bm{u}_i\|}{b_i} \right)^2 \quad \Longrightarrow \quad b_i \phi_i  \geqslant \|\bm{u}_i\|^2  ,  
\end{align}
Thus, the convex hull of \(S^i\) (see \cite{ceria1999convex} for proof) is given by:
\begin{align}
    \text{conv}(S^i) =&  \left\{ (\|\bm{u}_i\|, \phi_i, b_i) \in \mathbb{R}^3 : \ b_i \phi_i \geqslant \|\bm{u}_i\|^2, \right. \nonumber\\ 
    & \left. \ \ b_i \|\bm{u}_{\min}\| \ \leqslant \ \|\bm{u}_i\| \ \leqslant \  b_i \|\bm{u}_{\max}\|, \right. \nonumber \\ 
    & \left. \ \ 0 \leqslant b_i \leqslant 1,\quad  \phi_i \geqslant 0 \right\}.
\end{align}
\begin{figure}[h]
    \centering
    \includegraphics[width=0.9\linewidth]{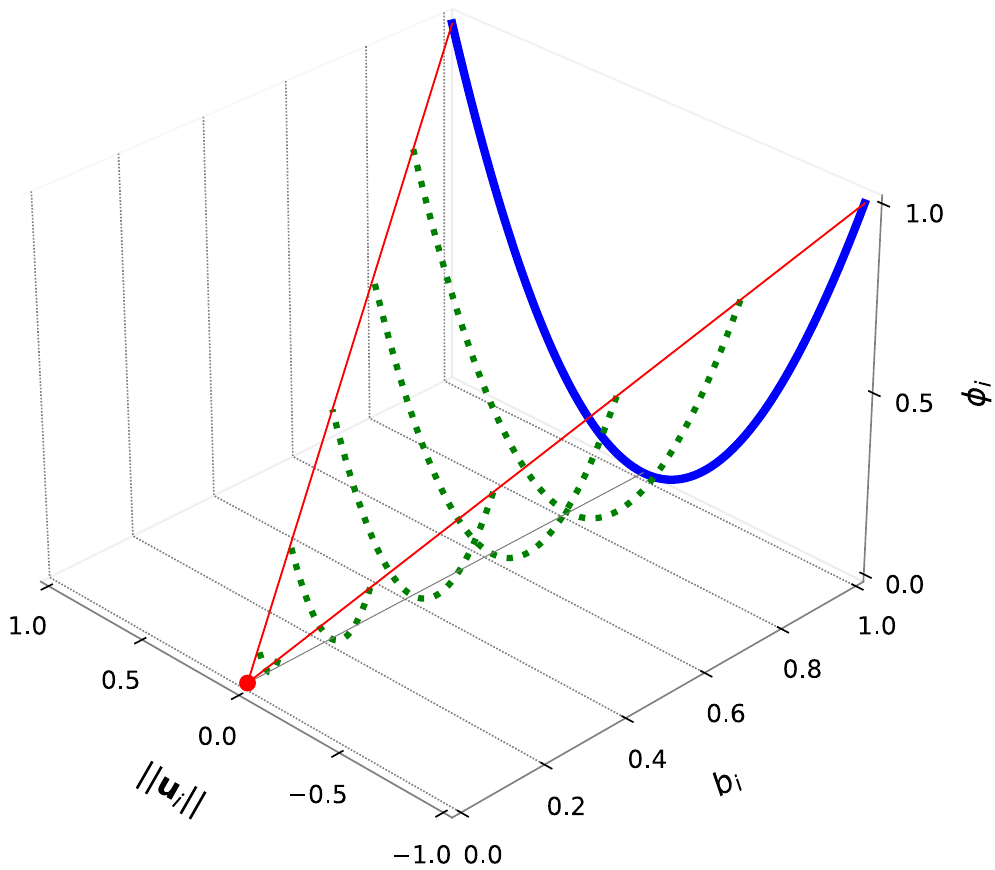}
    \caption{Feasible region of the perspective convex reformulation of a quadratic function.}
    \label{Fig:perspective}
\end{figure}
A graphical illustration of this convex hull for varying fractional values of \(b_i\), resulting from the perspective reformulation, is provided in Fig. \ref{Fig:perspective}. We apply this convex reformulation to tighten the \(\mathtt{RNLP}\) formulation.
\begin{flalign*}
    \mathtt{RNLP-T} ~~ \triangleq ~~ & \min ~~ \text{Eq. \eqref{eq:perspective-cuta}} ~~~ \text{ subject to: } \\
    & \text{Eq. \eqref{eq:RK-discrete}, \eqref{eq:RK-space}, \eqref{eq:RK-velocity}, \eqref{eq:inverse cubic distance}, \eqref{eq:initial conditions}, \eqref{eq:thruster}, \eqref{eq:proximity constraint}}, \eqref{eq:perspective-cutb} \\ 
    & \bm{r}_i, \bm{v}_i, \bm{u}_i \in \mathbb{R}^3, 0 \leqslant b_i \leqslant 1, \quad i = 0,\ldots,N
\end{flalign*}

\subsection{Rounding-based Feasible Solution}

Although the relaxed NLP formulation with perspective constraints provides a near-optimal and feasible solution, the integer variable solution still has infeasibility i.e. the variables have fractional values. To obtain a feasible integer solution, we fix the binary variables by rounding the fractional values from the perspective formulation solution to the nearest integer value as:
\begin{align}\label{eq:rounding}
    b_i = 
    \begin{cases}
        0, \quad b^*_i < \beta, \\
        1, \quad b^*_i \geqslant \beta 
    \end{cases}
\end{align}
\noindent We choose $\beta = 0.5$ as the default cut-off value. Other values for $\beta$ can also be used based on preferences for stricter cut-off values on thruster actions.

After fixing the binary variables, the optimization problem reduces to a NLP with only continuous variables ($\bm{r}_i,\bm{v}_i,\bm{u}_i$). We then solve the exact NLP to obtain a physically feasible integer solution to the space object tracking problem.
\begin{flalign*}
    \mathtt{NLP} ~~ \triangleq ~~ & \min ~~ \text{Eq. \eqref{eq:objective}} ~~~ \text{ subject to: } \\
    & \text{Eq. \eqref{eq:RK-discrete}, \eqref{eq:RK-space}, \eqref{eq:RK-velocity}, \eqref{eq:inverse cubic distance}, \eqref{eq:initial conditions}, \eqref{eq:thruster}, \eqref{eq:proximity constraint}}, \\
    & \bm{r}_i, \bm{v}_i, \bm{u}_i \in \mathbb{R}^3, \quad i = 0,\ldots,N
\end{flalign*}

The complete solution procedure for the space-tracking control and optimization problem is summarized in Algorithm \ref{alg:summary}.

\begin{algorithm}[ht]
\label{alg:summary}
  \caption{Space Object Tracking Optimization}\label{alg:SOTO}
  \KwIn{Primary Object Trajectory $\bm{\bar r_i}, \bm{\bar v_i}$}
  \KwOut{Secondary Object Thrusters $\bm{u_i}, b_i$}
\begin{enumerate}
    \item Solve the simulation \texttt{SIM} with zero thrust $\bm{u_i} = 0$
    \item Solve the relaxed problem with perspective constraint \texttt{RNLP-T} using the simulation solution as initial guess for the optimization
    \item Round and fix the relaxed binary variables (see eq. \eqref{eq:rounding}) returned by \texttt{RNLP-T} and re-solve the original continuous \texttt{NLP}, yielding an integer-feasible satellite-control schedule.
\end{enumerate}
\end{algorithm}

\section{NUMERICAL RESULTS}
\label{sec:numerical}
This section highlights the computational performance of the proposed problem formulation and algorithms for optimizing space object tracking via a detailed case study.

\subsection{Computational set-up and a case study}

In this case study, we implement the proposed algorithm to obtain an optimal trajectory for a space object tracking a primary object with a given predicted trajectory, i.e., position as a function of time. We solve a 1-hour time-horizon problem ($t_f - t_0 = 3600$ s) with ($\Delta t = 10$ s) intervals, resulting in $N=360$ intervals with 15,448 variables (15,088 continuous, 360 binary), 13,648 equality constraints, and 3,240 inequality constraints. The minimum ($\Delta_{\min}$) and maximum ($\Delta_{\max}$) distances for the secondary spacecraft are set at 10 km and 50 km, respectively. The initial distance between the space object and the secondary spacecraft is specified as $\bm{\Delta} = [10~ \text{km}, 10~ \text{km}, 10~ \text{km}]$. Realistic bounds on thruster outputs depend on the type of thrusters. Here, we use nominal values of $\bm{u}_{\min} = -1$ km/s$^2$ and $\bm{u}_{\max} = 1$ km/s$^2$ in each of the co-ordinates. The maximum number of time intervals during which the thrusters can be used is set at $N_b = 100$.

The optimization formulations and algorithms were modeled and implemented using JuMP v1.22 \cite{dunning2017jump} in Julia v1.10. Computational analyses utilized Ipopt 3.14.17 \cite{wachter2006implementation} as the continuous nonlinear optimization solver, running on an Apple M1 processor with 16GB of memory.

\subsection{Numerical Results}
\label{subsec:results}

The average total time to obtain a feasible ``optimal'' solution using our proposed strategy, including loading the primary object trajectory data and building the optimization models in JuMP, is approximately 38.2 sec. The simulation ($\mathtt{SIM}$) and exact formulation ($\mathtt{NLP}$) require an average of 0.6 and 3.6 seconds, respectively, while the relaxed binary optimization ($\mathtt{RNLP-T}$) converges in about 11.2 seconds.

\subsubsection{\underline{3D Optimal Trajectory}}

\begin{figure}[h]
    \centering
    \includegraphics[width=1.0\linewidth]{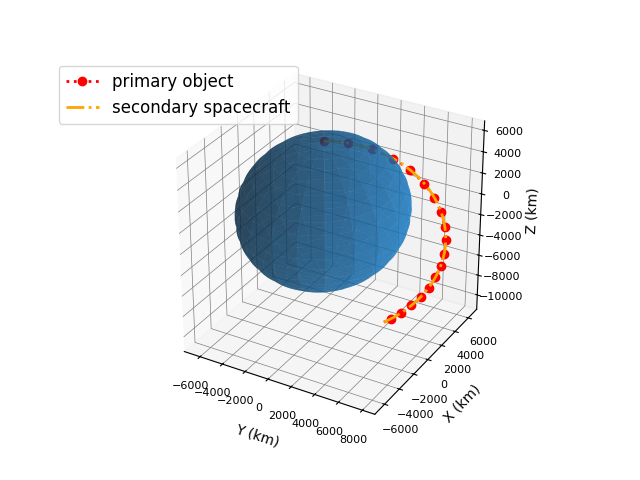}
    \caption{Optimal trajectory of the chaser and primary object around Earth.}
    \label{fig:3D-Trajectory}
\end{figure}

The optimal trajectory of the secondary spacecraft is sketched in Fig.~\ref{fig:3D-Trajectory} along with the pre-specified trajectory of the primary object in a Low-Earth Orbit.

As evident from the Fig.~\ref{fig:3D-Trajectory}, it is not possible to differentiate between the multiple trajectories, since the distance between the objects is very small compared to the size of the elliptical orbit of the primary.

\subsubsection{\underline{Inter-space Distance between Objects}}

\begin{figure}[h]
    \centering
    \includegraphics[width=0.9\linewidth]{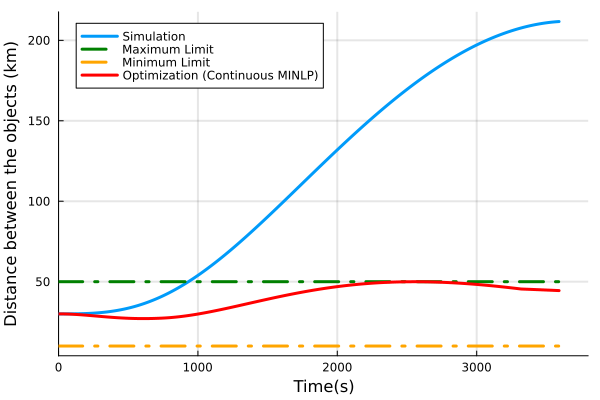}
    \caption{Distance (L1-norm) between the primary object and the secondary spacecraft.}
    \label{fig:Distance-plot}
\end{figure}

Fig.~\ref{fig:Distance-plot} presents the l1-norm distance between the primary object and the secondary spacecraft as a function of time for both the simulation and the integer feasible optimization ($\mathtt{NLP}$).

In the simulation, where thrusting is fixed to zero, the relative distance between the chaser spacecraft and the target object exceeds the 50 km upper limit at t = 1000s, reaching a maximum of 200 km — four times the upper bound. In contrast, in the optimization case with the same initial conditions, the secondary spacecraft remains within the 10 km to 50 km bounds from the primary object. This demonstrates the effectiveness of control actions from the optimization problem in satisfying the min/max proximity constraints.

\subsubsection{\underline{Perspective Constraint}}
\begin{figure}[h!]
    \centering
    \includegraphics[width=0.9\linewidth]{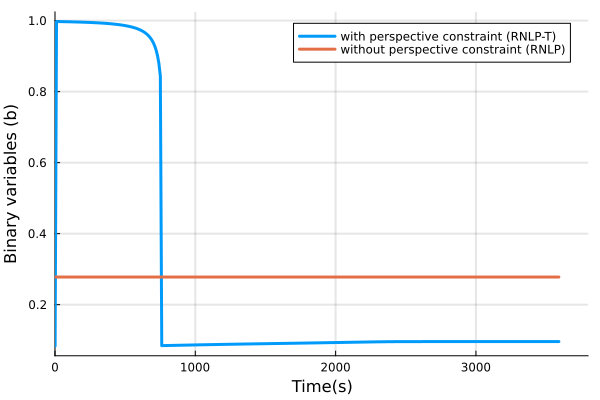}
    \caption{Binary solution of the relaxed NLP with and without perspective reformulation.}
    \label{fig:Binary-solution}
\end{figure}

Fig.~\ref{fig:Binary-solution} presents a comparison of the binary (on/off) solutions of the relaxed problems, both with ($\mathtt{RNLP-T}$) and without ($\mathtt{RNLP}$) the perspective reformulation from Sec. \ref{sec:perspective}.


The quality of a relaxed binary solution is determined by how close the binary variable values are to 0 or 1. The closer they are, the better the solution, as it indicates near feasibility and optimality. In Fig.~\ref{fig:Binary-solution}, the solution with the perspective constraint is of much higher quality compared to the solution without it. The relaxed binary variables in the perspective formulation are very close to either 1 ($b_i > 0.95$) or 0 ($b_i < 0.05$) making it much easier to find an integer feasible solution. The rounding algorithm applied to the perspective solution yields non-trivial, non-zero thrust profiles, while the solution without the perspective constraint ($b_i = 0.278$) leads to a trivial, impractical outcome, with all thrusters cut off.

\subsubsection{\underline{Thruster Profile}}
Applying the rounding criterion from Eq.~\ref{eq:rounding} to the relaxed binary solution with perspective reformulation (blue in Fig.~\ref{fig:Binary-solution}) results in a non-trivial thruster profile. The thrusters are switched on for the first 15 minutes (t=0 to t=900 s) and switched off from t=900 s to t=3600 s.

\begin{figure}[h]
    \centering
    \includegraphics[width=0.9\linewidth]{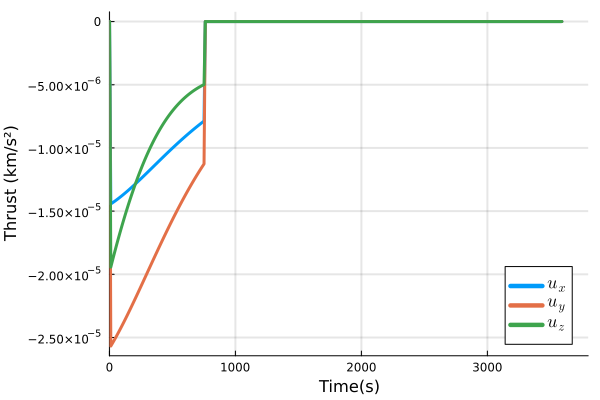}
    \caption{Optimal thruster profile in x, y and z directions.}
    \label{fig:thruster}
\end{figure}


The thruster action values, shown as continuous profiles in Fig.~\ref{fig:thruster}, minimize the operating cost and/or maximize the mission lifetime of the spacecraft tracking problem. The thruster values are orders of magnitude lower than other variables, such as coordinates and velocity. Notably, the highest thruster action is in the y-direction, with x and z-direction thruster values being comparable.

\subsection{Comments}
The optimization model in Section \ref{sec:form}, which includes discretized differential equations and operational constraints, is numerically unstable and ill-conditioned. The system is highly sensitive to parameters in the initial condition and variable or constraint bounds. To address these numerical issues, we use a simple reformulation by defining a new variable for the inverse cubic distance, as shown below:
\begin{align*}
    \Tilde{d}_{j,i} = \mu\,d_{j,i}
\end{align*}
Since a typical range for \(d\) is \(10^{-12} - 10^{-10}\) km\(^{-3}\), which is below the solver's tolerance of \(10^{-8}\), we replace it with \(\Tilde{d}\) in Eq.~\ref{eq:inverse cubic distance}, with values between \(10^{-6} - 10^{-4}\) s\(^{-2}\). 

Scaling the variables with nominal values for distance ($\bar r = 6378$ km i.e. earth radius) and acceleration ($\bar a = 10^{-3} $km s$^{-2}$) was implemented but didn't provide much improvement in numerical convergence or computational efficiency. Analysis on scaling the stiff nonlinear ODE for more robust and computationally efficient solution is part of ongoing work.

\section{CONCLUSIONS}
In this paper, we introduced a novel optimization approach for trajectory planning of a chaser spacecraft tracking a target object, while effectively handling collision avoidance, resource efficiency, and spatial constraints. By formulating the problem as a Mixed-Integer Nonlinear Program (MINLP) with nonlinear dynamics discretized using a fourth-order Runge-Kutta method, we eliminated the need for linearized approximations and extended the formulation to handle complex perturbations in the target's motion. The inclusion of continuous relaxation of discrete thrusting actions and perspective reformulation followed by rounding methods were found effective in solving this complex non-convex, non-linear problem efficiently within practical computation limits. The proposed method, validated through a numerical case study, demonstrates its suitability for onboard implementation in fast-moving objects, offering a promising solution for real-time spacecraft control in proximity operations. 

Future work could explore stochastic uncertainties (e.g., sensor noise), incorporate more sophisticated dynamics, and account for high-fidelity perturbations such as third-body gravity and atmospheric drag. Additionally, considering operational constraints like thruster latency and time-varying mass dynamics would enhance the approach’s robustness for a broader range of mission scenarios.

\addtolength{\textheight}{-12cm}   








\bibliographystyle{IEEEtran}
\bibliography{references}

\end{document}